\documentclass[12pt]{amsart}

\usepackage{amssymb,amscd}
\usepackage{amsbsy,amsmath,amsfonts,amsthm,extarrows,bm}
\usepackage[english]{babel}
\usepackage{graphicx}
\usepackage[all]{xy}
\usepackage[latin1]{inputenc}

\newtheorem{lemma}{Lemma}

\newtheorem{Atheorem}{Theorem}[section]
 
\newtheorem{Acorollary}[Atheorem]{Corollary}
 
\newenvironment{sproof}[1]
{\begin{proof}[#1]} {\end{proof}}

\newcommand{\Z}{\mathbb Z}

\newcommand{\sr}{\textnormal{sr}}
\newcommand{\GL}{\textnormal{GL}}
\newcommand{\SL}{\textnormal{SL}}
\newcommand{\D}{\textnormal{D}}
\newcommand{\El}{\textnormal{E}}
\newcommand{\rk}{\textnormal{rk}}
\newcommand{\U}{\textnormal{Um}}
\newcommand{\E}{\textnormal{E}}

\newcommand{\V}{\operatorname{V}}
\newcommand{\diag}{\operatorname{diag}}
\newcommand{\PID}{\operatorname{PID}}
\newcommand{\PIR}{\operatorname{PIR}}

\newcommand{\br}[1]{\lbrack #1 \rbrack}

\newcommand{\eb}{\mathbf{e}}
\newcommand{\gb}{\mathbf{g}}

\newcommand{\rb}{\mathbf{r}}

\newcommand{\mb}{\mathbf{m}}

\newcommand{\ZC}{ \mathbb{Z} \lbrack C \rbrack}
\newcommand{\ia}{\mathfrak{a}}
\newcommand{\ib}{\mathfrak{b}}
\newcommand{\ir}{\mathfrak{r}}

\title{Finitely generated modules over quasi-Euclidean rings}
\address{EPFL ENT CBS BBP/HBP. Campus Biotech. B1 Building, Chemin des mines, 9\\Geneva 1202, Switzerland}
\email{luc.guyot@epfl.ch}
\author{Luc Guyot}
\date{\today}
\thanks{The author thanks the Mathematics Institute of G\"oettingen and the Max Planck Institute for Mathematics in Bonn for the excellent conditions provided for his stay at these institutions, during which the paper was written.}

\begin{document}
\maketitle
\begin{abstract}
Let $R$ be a unital commutative ring and let $M$ be an $R$-module that is generated by $k$ elements but not less. Let $\El_n(R)$ be the subgroup of $\GL_n(R)$ generated by the elementary matrices. 
In this paper we study the action of $\E_n(R)$ by matrix multiplication on the set $\U_n(M)$ of unimodular rows of $M$ of length $n \ge k$.
Assuming $R$ is moreover Noetherian and quasi-Euclidean, e.g., $R$ is a direct sum of finitely many Euclidean rings,
we show that this action is transitive if $n > k$. 
We also prove that $\U_k(M) /\E_k(R)$ is equipotent with the unit group of $R/\ia_1$ where $\ia_1$ is the first invariant factor of $M$. 
These results encompass the well-known classification of Nielsen non-equivalent generating tuples in finitely generated Abelian groups.
\end{abstract}

\section{Introduction} \label{SecIntro}
In this paper rings are supposed unital and commutative. 
The unit group of a ring $R$ is denoted by $R^{\times}$. Let $M$ be a finitely generated $R$-module. 
We denote by $\rk_R(M)$ the minimal number of generators of $M$.
For $n \ge \rk_R(M)$, we denote by $\U_n(M)$ the set of \emph{unimodular rows} of $M$ of length $n$, i.e, the set of elements in $M^n$ whose components generate $M$. We consider the action of $\GL_n(R)$ on $\U_n(M)$ by matrix right-multiplication. Let $\E_n(R)$ be the subgroup of $\GL_n(R)$ generated by the \emph{elementary matrices}, i.e., the matrices which differ from the identity by a single off-diagonal element. 
Two unimodular rows $\mb$, $\mb' \in \U_n(M)$ are said to be \emph{$\El_n(R)$-equivalent} if there exists $E \in \El_n(R)$ such that 
$\mb' = \mb E$. Our chief concern is the description of the orbit set $\U_n(M)/\E_n(R)$ when $R$ enjoys a cancellation property 
shared by Euclidean rings, namely:
for every $n \ge 2$ and every $\rb = (r_1, \dots, r_n) \in R^n$, there exist $E \in \El_n(R)$ and $d \in R$ such that   
\begin{equation} \label{EqEnCancellation}
(d, 0,\dots,0) = \rb E.
\end{equation}
Rings with the above property are known as quasi-Euclidean rings in the sense of O'Meara and Cooke \cite{Ome65, Coo76} (see  \cite[Theorem 11]{AJLL14} for equivalence of definitions). 
A Noetherian quasi-Euclidean ring $R$ is therefore a \emph{principal ideal ring} (PIR), i.e., a ring whose ideals are principal. It is moreover an elementary divisor ring (see definition below) so that every finitely generated $R$-module admits an \emph{invariant factor decomposition}, that is a decomposition of the form 
\begin{equation} \label{EqInvariantFactorDecomposition}
R /\ia_1 \times R/\ia_2 \times  \cdots \times R /\ia_k
\text{ with } R \neq \ia_1 \supset \ia_2 \supset \cdots  \supset \ia_k
\end{equation}
where $k$ is necessarily equal to $\rk_R(M)$.
If $R$ is both Noetherian and quasi-Euclidean, properties (\ref{EqEnCancellation}) and (\ref{EqInvariantFactorDecomposition}) make the study of the action of $\E_n(R)$ on $\U_n(R)$ particularly amenable. Indeed, this assumption implies that $\El_n(R)$ acts transitively on $\U_n(M)$ for every $n > \rk_R(M)$ and enables us to exhibit a complete invariant of $\El_n(R)$-equivalence when $n = \rk_R(M)$. These two claims are corollaries of Theorem \ref{ThEnEquivalence} below.

\begin{Atheorem} \label{ThEnEquivalence}
Let $R$ be a Noetherian quasi-Euclidean ring and let $M$ be a finitely generated $R$-module.
Let
$
R/\ia_1  \times \cdots \times R/\ia_k
$
be the invariant factor decomposition of $M$.
Then every row in $\U_n(M)$ with $n \ge k$ is $\El_n(R)$-equivalent to a row of the form
$(\delta e_1, e_2, \dots, e_k, 0, \dots, 0)$ with $\delta \in (R/\ia_1)^{\times}$ and where $e_i \in M$ is defined by $(e_i)_i = 1 \in R/\ia_i$ and $(e_i)_j = 0$ for $j \neq i$. If $n > k$, then $\delta$ can be replaced by the identity.
\end{Atheorem}

\begin{Acorollary} \label{CorTransitivity}
Let $M$ be as in Theorem \ref{ThEnEquivalence}.
The action of $\E_n(R)$ on $\U_n(M)$ by matrix multiplication is transitive for every $n > \rk_R(M)$.
\end{Acorollary}

\begin{Acorollary} \label{CorDet}
Let $R$ be a Noetherian quasi-Euclidean ring.
Let $\ib_1, \dots, \ib_k$ be proper ideals of $R$ and set $\ib = \ib_1 + \cdots + \ib_k$.
Denote by $M$ the $R$-module $R/\ib_1  \times \cdots \times R/\ib_k$ and
for $\mb = (m_i)  \in M^{k}$ denote by $\det(\mb)$ the determinant of the matrix whose coefficients are the images in $R/\ib$ of the $(m_i)_j$'s 
via the natural maps $R/\ib_j \twoheadrightarrow R/\ib$.
Then $\mb, \mb' \in \U_k(M)$ are $\El_k(R)$-equivalent if and only if $\det(\mb) = \det(\mb')$.
\end{Acorollary}

Suppose $k = \rk_R(M)$. Then the ideal $\ib$ of Corollary \ref{CorDet} coincides with the first invariant factor $\ia_1$ of $M$ (cf. proof). 
Putting Theorem \ref{ThEnEquivalence} and Corollary \ref{CorDet} together, we see that $\mb \in \U_k(M)$ is $\El_k(R)$-equivalent to $(\det(\mb)e_1, e_2, \dots, e_k)$ where $(e_i)$ is as in Theorem \ref{ThEnEquivalence} and $\det$ as in Corollary \ref{CorDet}. Therefore the map $\gb \mapsto \det(\gb)$ induces a bijection from $\U_k(M) /\E_k(R)$ onto $(R/\ia_1)^{\times}$. The latter bijection endows $\U_k(M) /\E_k(R)$ with an Abelian group structure. The subject whether $\U_k(R) /\E_k(R)$ has a group structure for $R$ a commutative ring of finite Krull dimension is well studied 
\cite{SV76, Van83, Van89, Rao98, Fas11} but we don't know of any similar results for modules. By analogy with \cite[Definition 11.3.9]{McCR87}, we can define the elementary rank of a finitely generated $R$-module $M$, 
for any associative ring $R$ with identity, as the least integer $e$ such that the action of $\E_n(R)$ on $\U_n(M)$ is transitive for all $n > e$. This rank is not less than $\rk_R(M) - 1$ and not greater than $\rk_R(M) -1 + \sr(M)$, where $\sr(M)$ is the stable rank of $M$, a natural generalization of the Bass stable rank of rings to modules \cite[Definition 6.7.2]{McCR87}. We showed in our situation that the elementary rank is $\rk_R(M) - 1$ if $R/\ia_1$ has $2$ elements and coincides with $\rk_R(M)$ otherwise.

Specifying the above results to $R = \Z$ yields the characterization of Nielsen equivalent generating tuples in finitely generated Abelian groups. This characterization was obtained in part by several authors \cite{NN51, LM93, DG99} and reaches its complete form in \cite{Oan11}.
In order to present it, we introduce the following definitions.
Given a finitely generated group $G$, denote by $\rk(G)$ the minimal number of generators of $G$. For $n \ge \rk(G)$, let $\V_n(G)$ be the set of \emph{generating $n$-vectors} of $G$, i.e., the set of elements in $G^n$ whose components generate $G$.
Two generating  $n$-vectors are said to be \emph{Nielsen equivalent} if they can be related by a finite sequence of  transformations of $G^n$ taken in the set $\{ L_{ij}, I_i; 1 \le i \neq j \le n\}$ where
$L_{ij}$ and $I_i$ replace the component $g_i$ of $\gb = (g_1, \dots, g_n) \in G^n$ by $g_j g_i$ and $g_i^{-1}$ respectively and leave the other components unchanged.

\begin{Acorollary} \label{CorNielsenAbel}
Let $G$ be a finitely generated Abelian group whose invariant factor decomposition is
$$
\Z_{d_1} \times \cdots  \times \Z_{d_k}
$$
with $1 \neq d_1 \, \vert \, d_2 \, \vert \, \cdots \, \vert \, d_k, \, d_i \ge 0$ and where $\Z_{d_i}$ stands for $\Z/d_i\Z$ 
(in particular $\Z_0 = \Z$).
Then every generating $n$-vector $\gb$ with $n \ge k$ is
Nielsen equivalent to
$(\delta e_1, e_2, \dots, e_k, 0, \dots, 0)$ for some $\delta \in (\Z_{d_1})^{\times}$ 
and where $e_i \in G$ is defined by $(e_i)_i = 1 \in \Z_{d_i}$ and $(e_i)_j = 0$ for $j \neq i$. 

\begin{itemize}
\item If $n > k$, then $\delta$ above can be replaced by the identity.
\item If $n = k$ then $\delta$ must be $\pm \det(\gb)$ with $\det$ defined as in Corollary \ref{CorDet}.
\end{itemize}

In particular $G$ has only one Nielsen equivalence class of generating $n$-vectors for $n > k$ while it has $\max(\varphi(d_1)/2, 1)$ Nielsen equivalence classes of generating $k$-vectors where $\varphi$ denotes the Euler totient function extended by $\varphi(0) = 0$.
\end{Acorollary}

Our results allow further applications to the study of Nielsen equivalence in split extensions of Abelian groups by cyclic or free Abelian groups. Consider for instance an infinite cyclic group $C$ and denote by $\ZC$ its integral group ring.
Let $R$ be a quasi-Euclidean quotient of $\ZC$, e.g., $R = \Z_n\br{C}$ for $n$ a square-free integer, and let $M$ be a finitely generated $R$-module. Then the image of $C$ in $R$ is a subgroup of $R^{\times}$ so that $C$ acts naturally on $M$ by automorphisms.
Let $G = M \rtimes C$ be the corresponding semi-direct product.
 Let $T$ be the subgroup of $R^{\times}$ generated by the images of $-1$ and $C$, set $\Lambda = R/\ia_1$ where $\ia_1$ is the first invariant factor of $M$, and let $T_{\Lambda}$ be the image of $T$ in $\Lambda^{\times}$.
Our results allow us to show that the set of Nielsen equivalence classes of generating $k$-tuples of $G$ is equipotent with $\Lambda^{\times}/T_{\Lambda}$ for $k = \rk(G)$ \cite{Guy16b}. (See also \cite[Corollary 4.$ii$]{Guy16a} for an application of Corollary \ref{CorTransitivity}).

We introduce now definitions and preliminary results that will be used by the proofs of our statements in Section \ref{SecProofs}.

A ring $R$ is said to be an \emph{elementary divisor ring} if every matrix $A$ over $R$ admits a \emph{diagonal reduction}, i.e., if we can find invertible matrices $P,\,Q$ and elements $d_1,\dots,d_n \in R$ such that $PAQ = \diag(d_1,\dots,d_n)$ and $d_1 \, \vert \, d_2 \, \vert \cdots \vert \, d_n$. Principal ideal domains (PID) are elementary divisor rings \cite[Theorem 4]{DuFo04}. This classical result extends effortlessly to PIRs thanks to a theorem of Hungerford. Indeed, the class of elementary divisor rings is stable under taking quotients and direct sums. As a $\PIR$ is a direct sum of rings, each of which is the homomorphic image of a $\PID$ \cite[Theorem 1]{Hun68}, our claim follows. 
A Noetherian elementary divisor ring has moreover a unique invariant factor decomposition \cite[Theorems 9.1 and 9.3]{Kap49}, thus we showed
\begin{lemma} \label{LemPIRIsEDR}
Let $R$ be a $\PIR$. Then the following hold:
\begin{itemize}
\item[$(i)$] $R$ is an elementary divisor ring.
\item[$(ii)$] Every finitely generated $R$-module $M$ has a unique invariant factor decomposition, i.e., a decomposition of the form
$R /\ia_1 \times R/\ia_2 \times  \cdots \times R /\ia_n$ with $R \neq \ia_1 \supset \ia_2 \supset \cdots  \supset \ia_n$ where
the factors $\ia_i$ are uniquely determined by the latter condition. For such decomposition, $n$ is the minimal number of generators of $M$, that is $n = \rk_R(M)$.
\end{itemize}
\end{lemma}

The following result of Whitehead \cite[Example I.1.11]{Wei13} will come in handy in the proof of Theorem \ref{ThEnEquivalence}.
Let $R$ be a unital commutative ring. For $u \in R^{\times}$, we denote by $D_i(u) \in \GL_n(R)$ the diagonal matrix which coincides with the identity except possibly for its $(i,i)$-entry, which is set to $u$. Let $\D_n(R^{\times})$ be the subgroup of $\GL_n(R)$ generated by the matrices $D_i(u)$. Whitehead's lemma implies that a matrix $D \in \D_n(R^{\times})$ lies in $\SL_n(R)$ if and only if 
 it lies in $\El_n(R)$.

Corollary \ref{CorDet} elaborates on Theorem \ref{ThEnEquivalence} by showing that the unit $\delta$ of the theorem identifies with a natural invariant of $\El_n(R)$-equivalence, namely the determinant in the largest quotient of $M$ that is a free module. 
This invariant extends Diaconis-Graham's invariant defined for finitely generated Abelian groups \cite{DG99}.
It can be defined for any commutative unital ring $R$ and any finitely generated $R$-module $M$. Consider a generating set 
$\mb = (m_1, \dots , m_n)$ of $M$ with minimal cardinality. We say that $r \in R$ is \emph{involved} in a relation of $M$ with respect to $\mb$ if there is $(r_i) \in R^n$ such that $\sum r_i m_i = 0$ and $r = r_i$ for some $i$. 
Denote by $\ir(M)$ the set of elements of $R$ which are involved in a relation of $M$ with respect to $\mb$. 
Clearly, $\ir(M)$ is an ideal of $R$ and it is easily checked that $\ir(M)$ is independent of $\mb$. 
Let $\overline{\mb} = \pi(\mb)$ be the image of $\mb$ by the natural map 
$\pi: M \twoheadrightarrow M/\ir(M) M$ and let $\eb$ be the canonical basis of $(R/\ir(M))^n$. 
Then the map $\overline{\mb} \mapsto \eb$ induces an isomorphism 
$\varphi_{\mb}$ from $M/\ir(M) M$ onto $(R/\ir(M))^n$. For $\mb' \in \U_n(M)$, 
we define $\det_{\mb}(\mb')$ as the determinant of 
$\varphi_{\mb} \circ \pi(\mb')$ with respect to $\eb$. Let $\varphi \Doteq \varphi_{\mb} \circ \varphi^{-1}_{\mb'}$. 
Because the identity $\det_{\mb} = \det(\varphi) \det_{\mb'}$ holds, we have shown the following
\begin{lemma} \label{LemUnitRelated}
Let $n = \rk_R(M)$ and let $\mb, \mb' \in \U_n( M)$. There exists $u \in (R/\ir(M))^{\times}$ such that $\det_{\mb} = u \det_{\mb'}$.
\end{lemma}

It is straightforward to check that $\det_{\mb}$ is an invariant of $\El_n(R)$-equivalence, i.e., 
$\det_{\mb}(\mb'E) = \det_{\mb}(\mb')$ for every $E \in \El_n( R)$.
If $\ir(M) = R$, then $\det_{\mb}$ is trivial and hence useless. 
This doesn't happen when $R$ is an elementary divisor ring, for $\ir(M)$ is then the first invariant factor $\ia_1$ of $M$. The identity elements of each factor ring in a decomposition $M \simeq R/\ia_1  \times \cdots \times R/\ia_n$ where $n = \rk_R(M)$ form a unimodular row of $M$. We refer to this unimodular row as the \emph{unimodular row naturally associated}  to the given decomposition.\\

\paragraph{\textbf{Acknowledgements}.}
The author is grateful to Pierre de la Harpe, Tatiana Smirnova-Nagnibeda and Laurent Bartholdi for their encouragements and useful comments made on preliminary versions of this paper. He is also thankful to Wolfgang Pitsch and Daniel Oancea for their early corrections and helpful suggestions.

\section{Proofs} \label{SecProofs}
\begin{sproof}{Proof of Theorem \ref{ThEnEquivalence}}
The case $k = 1$ follows straighforwardly from (\ref{EqEnCancellation}). 
We assume now that $k \ge 2$ and reverse the order of the sequence $(\ia_i)$ for notational convenience, 
supposing that $\ia_j \subset \ia_{j + 1}$ for every $j$. 
Let $\mb = (m_i) \in \U_n(M)$ with $n \ge k$. 
We set $m_{ij} \Doteq (m_i)_j$ for every $1 \le i \le n,\, 1 \le j \le k$ and identify $\mb$ with the matrix $(m_{ij})$. 
Applying (\ref{EqEnCancellation}) 
to every row $\mb_i \Doteq (m_{ij})_{i \le j \le k}$ we obtain a matrix $E \in \El_n(R)$ 
such that $\mb' = \mb E$ is a lower-triangular matrix $(m'_{ij})$ with $m'_{11} = 1$. 
We claim that $m'_{jj}$ is a unit of $R/\ia_j$ for every $2 \le j \le k$. 
To see this, we consider the image $\rho_l$ of $(m'_{ij})_{1 \le i,\, j \le l}$ via the natural map induced 
by the projection $R/\ia_{l - 1} \twoheadrightarrow R/\ia_l$  
for every $2 \le l \le k$. Since the rows of 
$\rho_l$ generate $(R/\ia_l)^l$, the matrix $\rho_l$ is invertible \cite[Theorem 2.4]{Mat89}.
Therefore $\det(\rho_l) \in (R/\ia_l)^{\times}$ and subsequently $m'_{ll} \in (R/\ia_l)^{\times}$.
 As a result we readily find $E' \in \El_n(R)$ 
such that $\mb'' \Doteq \mb'E'= \diag(1,\, m''_2,\, \dots,\,m''_k)$ with $m''_j \Doteq m'_{jj}$ for $2 \le j \le k$. 
If $k = 2$, we are done.
Otherwise we consider the matrices $E_j \Doteq D_j(1/m''_j)D_k(m''_j)$ for $2 \le j \le k - 1$. Since $E_j \in \El_k(R/\ia_j)$ for every $j$ by Whitehead's lemma, each matrix $E_j$ has a lift in $E_k(R)$ and the product of these lifts
is a matrix $E'' \in \El_n(R)$ satisfying $\mb''E'' = \diag(1, \dots, 1,  m'''_k)$ with $m'''_k \in (R/\ia_k)^{\times}$. 
If $n > k$, we can store $1 - m'''_k$ in the $(k, k + 1)$-entry of $\mb''$, then turn $m'''_k$ into $1$ 
and eventually cancel the $(k, k + 1)$-entry with obvious elementary column transformations.
\end{sproof}

\begin{sproof}{Proof of Corollary \ref{CorDet}}
Since $\det$ is invariant under elementary row transformations, the 'only if' part is established. To prove the converse, it suffices to show that there is $u \in (R/\ib)^{\times}$ such that for every $\mb \in M^{k}$,  $\det(\mb) = u \delta(\mb)$ where $\delta = \delta(\mb)$ is the unit given by Theorem \ref{ThEnEquivalence}. Let $\pi,\,\pi':M \twoheadrightarrow (R/\ib)^k$ be the $R$-epimorphisms naturally induced by the invariant factor decomposition and the decomposition $R/\ib_1  \times \cdots \times R/\ib_k$ respectively. Our claim certainly holds if  $\pi' = \varphi \circ \pi$ for some automorphism $\varphi$ of $(R/\ib)^k$.

To prove the latter fact, we first show that $\ib$ is the first invariant factor in the decomposition of $M$ 
if $k = \rk_R(M)$ and $\ib = R$ if $k > \rk_R(M)$. 
As $R$ is a $\PIR$, we can write $\ib_i = b_i R$ with $b_i \in R$ for every $i$. 
Let $A = \diag(b_1, \dots, b_k)$, with $A$ square of order $k$. 
Since $R$ is an elementary divisor ring, the matrix $A$ 
admits a diagonal reduction $\diag(d_1, \dots, d_k)$. As the ideal generated by 
the coefficients of $A$ is invariant under this reduction, 
we have $\ib = d_1 R$. If $d_1 \notin R^{\times}$, the ideals $d_i R$ correspond 
to the invariant factors $\ia_i$ of $M$ and hence $k = \rk_R(M)$. 
Otherwise $\ib = R$ and $\det$ vanishes on $\U_k(M)$ while rows in 
$\U_k(M)$ are all $\E_k(R)$-equivalent by Theorem \ref{ThEnEquivalence}. 
Therefore we can assume that $k = \rk_R(M)$ and $\ib = d_1 R$ is a proper ideal of $R$.

Considering the unimodular rows naturally associated to our two decompositions of $M$ as a direct sum of cyclic factors, 
we see that the existence of $\varphi$ is established by Lemma \ref{LemUnitRelated}. 
This proves the claim and hence the result.
\end{sproof}

\begin{sproof}{Proof of Corollary \ref{CorNielsenAbel}}
The group $G$ is a $\Z$-module and $\U_n(G)$ naturally identifies with $\V_n(G)$. Considering the matrix counterparts of the transformations $L_{ij}$ and $I_i$, it is easily checked that the Nielsen classes in $\V_n(G)$ concide with the orbits in $\U_n(G)$ of 
$\GL_n(\Z) = \D_n(\{\pm 1\}) \El_n(\Z)$ where $\D_n(\{\pm 1\})$ is the group of diagonal matrices with diagonal entries in $\{ \pm 1\}$. The result is then a straightforward consequence of 
Theorem \ref{ThEnEquivalence} and Corollary \ref{CorDet}.
\end{sproof}

\bibliographystyle{alpha}
\bibliography{Biblio}

\def\cprime{$'$} \def\cprime{$'$} \def\cprime{$'$} \def\cprime{$'$}
  \def\cprime{$'$}
\begin{thebibliography}{AJLL14}

\bibitem[AJLL14]{AJLL14}
Adel Alahmadi, S.~K. Jain, T.~Y. Lam, and A.~Leroy.
\newblock Euclidean pairs and quasi-{E}uclidean rings.
\newblock {\em J. Algebra}, 406:154--170, 2014.

\bibitem[Coo76]{Coo76}
George~E. Cooke.
\newblock A weakening of the {E}uclidean property for integral domains and
  applications to algebraic number theory. {I}.
\newblock {\em J. Reine Angew. Math.}, 282:133--156, 1976.

\bibitem[DF04]{DuFo04}
David~S. Dummit and Richard~M. Foote.
\newblock {\em Abstract algebra}.
\newblock John Wiley \& Sons, Inc., Hoboken, NJ, third edition, 2004.

\bibitem[DG99]{DG99}
Persi Diaconis and Ronald Graham.
\newblock The graph of generating sets of an abelian group.
\newblock {\em Colloq. Math.}, 80(1):31--38, 1999.

\bibitem[Fas11]{Fas11}
Jean Fasel.
\newblock Some remarks on orbit sets of unimodular rows.
\newblock {\em Comment. Math. Helv.}, 86(1):13--39, 2011.

\bibitem[Guy16a]{Guy16a}
L.~Guyot.
\newblock Generators in split extensions of {A}belian groups by cyclic groups.
\newblock Preprint, arXiv:1604.08896 [math.GR], 2016.

\bibitem[Guy16b]{Guy16b}
L.~Guyot.
\newblock Generators of split metabelian groups.
\newblock In preparation, 2016.

\bibitem[Hun68]{Hun68}
Thomas~W. Hungerford.
\newblock On the structure of principal ideal rings.
\newblock {\em Pacific J. Math.}, 25:543--547, 1968.

\bibitem[Kap49]{Kap49}
Irving Kaplansky.
\newblock Elementary divisors and modules.
\newblock {\em Trans. Amer. Math. Soc.}, 66:464--491, 1949.

\bibitem[LM93]{LM93}
Martin Lustig and Yoav Moriah.
\newblock Generating systems of groups and {R}eidemeister-{W}hitehead torsion.
\newblock {\em J. Algebra}, 157(1):170--198, 1993.

\bibitem[Mat89]{Mat89}
Hideyuki Matsumura.
\newblock {\em Commutative ring theory}, volume~8 of {\em Cambridge Studies in
  Advanced Mathematics}.
\newblock Cambridge University Press, Cambridge, second edition, 1989.
\newblock Translated from the Japanese by M. Reid.

\bibitem[MR87]{McCR87}
J.~C. McConnell and J.~C. Robson.
\newblock {\em Noncommutative {N}oetherian rings}.
\newblock Pure and Applied Mathematics (New York). John Wiley \& Sons, Ltd.,
  Chichester, 1987.
\newblock With the cooperation of L. W. Small, A Wiley-Interscience
  Publication.

\bibitem[NN51]{NN51}
Bernhard~H. Neumann and Hanna Neumann.
\newblock Zwei {K}lassen charakteristischer {U}ntergruppen und ihre
  {F}aktorgruppen.
\newblock {\em Math. Nachr.}, 4:106--125, 1951.

\bibitem[Oan11]{Oan11}
Daniel Oancea.
\newblock A note on {N}ielsen equivalence in finitely generated abelian groups.
\newblock {\em Bull. Aust. Math. Soc.}, 84(1):127--136, 2011.

\bibitem[O'M65]{Ome65}
O.~T. O'Meara.
\newblock On the finite generation of linear groups over {H}asse domains.
\newblock {\em J. Reine Angew. Math.}, 217:79--108, 1965.

\bibitem[Rao98]{Rao98}
Ravi~A. Rao.
\newblock An abelian group structure on orbits of ``unimodular squares'' in
  dimension {$3$}.
\newblock {\em J. Algebra}, 210(1):216--224, 1998.

\bibitem[vdK83]{Van83}
Wilberd van~der Kallen.
\newblock A group structure on certain orbit sets of unimodular rows.
\newblock {\em J. Algebra}, 82(2):363--397, 1983.

\bibitem[vdK89]{Van89}
Wilberd van~der Kallen.
\newblock A module structure on certain orbit sets of unimodular rows.
\newblock {\em J. Pure Appl. Algebra}, 57(3):281--316, 1989.

\bibitem[VS76]{SV76}
L.~N. Vaser{\v{s}}te{\u\i}n and A.~A. Suslin.
\newblock Serre's problem on projective modules over polynomial rings, and
  algebraic {$K$}-theory.
\newblock {\em Izv. Akad. Nauk SSSR Ser. Mat.}, 40(5):993--1054, 1199, 1976.

\bibitem[Wei13]{Wei13}
Charles~A. Weibel.
\newblock {\em The {$K$}-book}, volume 145 of {\em Graduate Studies in
  Mathematics}.
\newblock American Mathematical Society, Providence, RI, 2013.
\newblock An introduction to algebraic $K$-theory.

\end{thebibliography}

\end{document}